\documentclass[letterpaper]{article}

\usepackage{graphicx} 
\usepackage[section]{placeins} 

\usepackage[aboveskip=0pt,belowskip=-2pt]{subcaption} 
\usepackage[]{bm} 
\usepackage[]{bbold} 
\usepackage[]{amsmath} 

\usepackage[breaklinks, colorlinks=false,linkcolor=black, citecolor=blue, urlcolor=blue]{hyperref}
\usepackage{authblk}
\usepackage[capitalise]{cleveref} 

\usepackage[]{footmisc}

\title{A multiresolution adaptive wavelet method for nonlinear partial differential equations}
\author{Cale Harnish$^{1}$, Luke Dalessandro$^{2}$, Karel
  Matou\v{s}$^{1}$\footnote{Corresponding author.}\; and Daniel Livescu$^{3}$}

\affil{
${}^{1}$ University of Notre Dame, Notre Dame, IN, USA\\
${}^{2}$ Indiana University, Bloomington, IN, USA\\
${}^{3}$ Los Alamos National Laboratory, Los Alamos, NM, USA\\
${}$\\
charnish@nd.edu, ldalessa@iu.edu, kmatous@nd.edu, livescu@lanl.gov\\
}
\begin{document}
\maketitle
\begin{abstract}
The multiscale complexity of modern problems in computational science
and engineering can prohibit the use of traditional numerical methods
in multi-dimensional simulations. Therefore, novel algorithms are
required in these situations to solve partial differential equations
(PDEs) with features evolving on a wide range of spatial and temporal
scales. To meet these challenges, we present a multiresolution wavelet
algorithm to solve PDEs with significant data compression and explicit
error control. We discretize in space by projecting fields and spatial
derivative operators onto wavelet basis functions. We provide error
estimates for the wavelet representation of fields and their
derivatives. Then, our estimates are used to construct a sparse
multiresolution discretization which guarantees the prescribed
accuracy. Additionally, we embed a predictor-corrector procedure
within the temporal integration to dynamically adapt the computational
grid and maintain the accuracy of the solution of the PDE as it
evolves. We present examples to highlight the accuracy and adaptivity
of our approach.
\end{abstract}

\section{Introduction}
Modern computational science and engineering applications are inherently multiphysics and multiscale. For example, models of the global ocean \cite{ocean}, detonation combustion \cite{detonation}, asteroid impacts \cite{asteroid}, mechanics of materials \cite{matousJCP}, and supernova remnants \cite{supernova} all must solve partial differential equations (PDEs) with spatial and temporal scales across many orders of magnitude. Several novel numerical methods have been developed to address this computational challenge. For example, adaptive mesh refinement (AMR) \cite{Berger1984AdaptiveEquations, Fatkullin2001AdaptiveProblems}, multigrid methods \cite{Brandt1977Multi-LevelProblems, Hackbusch1978OnEquations,matousDewen}, Chimera overset grids \cite{chimera}, and remeshing/refining finite element methods (FEM) \cite{Dong2003P-refinementP-threads, Gui1986TheDimension, Gui1986TheDimensionb, Rajagopal2007AInterfaces} have been used to accomplish a great deal of contemporary computational modeling. However all of these methods become computationally expensive when the user does not know \emph{a priori} the spatial and temporal locations of interesting solution features. In this work, we propose a wavelet based method which is well-suited for problems with dynamic spatial and temporal scales.
 
Wavelet based numerical methods have been shown to be efficient for modeling multiscale and multiphysics problems because they provide spatial adaptivity through the use of multiresolution basis functions \cite{Jawerth1994AnAnalyses, Schneider2010WaveletDynamics}. Furthermore, current wavelet solvers have achieved several notable accomplishments, including: significant data compression \cite{Liandrat1990ResolutionApproximation, Beylkin1997OnBases, Bertoluzza1996AdaptiveEquation}, bounded energy conservation \cite{Ueno2003AProperty, Qian1993WaveletsEquations}, modeling stochastic systems \cite{Kong2016NonlinearTechnique}, multiscale model reduction \cite{RodyEtAl}, and solving coupled systems of nonlinear PDEs \cite{Paolucci2014WAMR:Algorithm, Paolucci2014WAMR:Algorithmb, Nejadmalayeri2015ParallelPDEs, Dubos2013AGrids, Sakurai2017CoherentWavelets}. However, some implementations only solve PDEs in infinite or periodic domains (\emph{e.g.}, \cite{Frohlich1994AnComputations, Goedecker2009WaveletsPhysics, Iqbal2014AnModeling}), some do not utilize the data compression ability of wavelets, resulting in a costly uniform grid (\emph{e.g.}, \cite{Qian1993WaveletsEquations, Le2015Reduced-orderWinds, Lin2001ConnectionEquation}), and some use finite difference methods to calculate the spatial derivatives, inhibiting the ability to solve PDEs in the wavelet domain and control accuracy (\emph{e.g.}, \cite{Paolucci2014WAMR:Algorithm, Paolucci2014WAMR:Algorithmb, Nejadmalayeri2015ParallelPDEs, Holmstrom1999SolvingWavelets}). 
 
To overcome some of the difficulties mentioned above, we have developed an algorithm which retains the advantages of other wavelet methods while attempting to overcome their limitations. Specifically, this work extends our one-dimensional solver described in \cite{Harnish2018AdaptiveControl} into multiple spatial dimensions. Since our proposed numerical method exploits the properties of wavelet basis functions, it is helpful to provide a brief outline of wavelet principles. Therefore in \cref{sec:theory}, we summarize the creation of wavelet basis functions and define the operations needed to solve PDEs using this basis. Then, in \cref{sec:algorithm} we describe the numerical implementation and in \cref{sec:verification_2D} we present illustrative numerical examples.
\section{Wavelet Theory}
\label{sec:theory}
A multiresolution analysis (MRA) provides the formal mathematical framework for a wavelet family of basis functions \cite{Daubechies1992TenDaubechies}. An MRA of a domain $\Omega$ consists of a progression of nested approximation spaces $V_{j}$ and their associated dual spaces $\widetilde{V}_{j}$ such that the union of these spaces is the $L^{2}(\Omega)$ space \cite{Cohen2000MultiscaleDomains}. The wavelet spaces $W_{j}$, and their associated dual spaces $\widetilde{W}_{j}$ are then defined as the complements of the approximation spaces $V_{j}$ in $V_{j+1}$ \cite{Qian1993WaveletsEquations,Bacry1992AEquations},
\begin{align} 
\label{eqn:MRA_spaces}
    V_{j} &\subset V_{j+1}, &  \overline{\bigcup_{j} V_{j}} =& L^{2}(\Omega), & V_{j+1} &= V_{j} \oplus W_{j}.
\end{align}
Then, multidimensional representations are defined by tensor products. For example, the two-dimensional space $\bm{V}_{j}$ is defined by the tensor product of two one-dimensional spaces $V_{j}$:
\begin{align}
\label{eqn:MRA_2D}
    \bm{V}_{j+1} &= V_{j+1} \otimes V_{j+1}, \nonumber\\
    \bm{V}_{j+1} &= \left( V_{j} \oplus W_{j} \right) \otimes\left( V_{j} \oplus W_{j} \right), \nonumber\\
    \bm{V}_{j+1} &= \underbrace{\left( V_{j} \otimes V_{j} \right)}_{\lambda = 0} \oplus \underbrace{\left( W_{j} \otimes V_{j} \right)}_{\lambda = 1} \oplus \underbrace{\left( V_{j} \otimes W_{j} \right)}_{\lambda = 2} \oplus \underbrace{\left( W_{j} \otimes W_{j} \right)}_{\lambda = 3}.
\end{align}
Therefore, the MRA creates four types of two-dimensional basis where each is designated by $\lambda$ and defined by the appropriate products of one-dimensional basis (\emph{i.e.}, $\phi_{k}^{j} (x)$ and $\psi_{k}^{j} (x)$) \cite{Daubechies1992TenDaubechies}. Note that the multiresolution nature of wavelets requires the use of two types of indices. One to define the resolution level $j$, and another to define the spatial locations $k$ on a particular resolution level. 

In general, wavelet bases do not have a closed-form expression, instead they are defined in terms of four types of filter coefficients (\emph{i.e.}, $h_{i}, \ \widetilde{h}_{i}, \ g_{i},$ and $\widetilde{g}_{i}$) \cite{Goedecker2009WaveletsPhysics, deVilliers2003Dubuc--DeslauriersInterval}. Our algorithm uses the Deslauriers-Dubuc wavelet family, with second generation wavelets near spatial boundaries, as defined in \cite{deVilliers2003Dubuc--DeslauriersInterval}. Furthermore, a single parameter $p$ defines the properties of this basis, such as the number of vanishing moments and the degree of continuity \cite{Harnish2018AdaptiveControl}. We discretize space by projecting each continuous field $f(\vec{x})$ onto the wavelet basis $\bm{\phi}_{\vec{k}}^{0} (\vec{x})$ and ${}^{\lambda}\bm{\psi}_{\vec{k}}^{j} (\vec{x})$, where $\vec{\bullet}$ indicates a vector. The corresponding wavelet coefficients $s_{\vec{k}}^{0}$ and ${}^{\lambda}d_{\vec{k}}^{j}$ are defined by integrating the field with the dual basis,
\begin{align}
\label{eqn:coeff_def}
    s_{\vec{k}}^{0} = \int_{\Omega} f(\vec{x}) \ \bm{\widetilde{\phi}}_{\vec{k}}^{0}(\vec{x}) \ \mathrm{d} \Omega \qquad \mathrm{and} \qquad {}^{\lambda}d_{\vec{k}}^{j} = \int_{\Omega} f(\vec{x}) \ {}^{\lambda}\bm{\widetilde{\psi}}_{\vec{k}}^{j}(\vec{x}) \ \mathrm{d} \Omega.
\end{align}
Leveraging properties of the Deslauriers-Dubuc wavelet family, the integrals in \cref{eqn:coeff_def} can be solved exactly and are replaced with the matrix operator $\bm{F}$, defined in terms of the filter coefficients, as shown in \cite{Harnish2018AdaptiveControl}. Repeated application of this operator yields all of the wavelet coefficients on each resolution level. As the resolution level increases, the magnitude of the calculated ${}^{\lambda}d_{\vec{k}}^{j}$ coefficients decreases, and on some resolution level $j = j_{\mathrm{max}}$, all of the coefficients will be below a prescribed tolerance $\varepsilon$. It has been shown by many authors (\emph{e.g.}, \cite{Paolucci2014WAMR:Algorithm, Nejadmalayeri2015ParallelPDEs, Holmstrom1999SolvingWavelets}) that discarding those ${}^{\lambda}d_{\vec{k}}^{j}$ coefficients with a magnitude less than $\varepsilon$, results in the discretization
\begin{align}
\label{eqn:multi}
    f_{\varepsilon}(\vec{x}) = \sum_{\vec{k}} s_{\vec{k}}^{0} \ \bm{\phi}_{\vec{k}}^{0}(\vec{x}) + \sum_{j=1}^{j_{\mathrm{max}}} \sum_{\lambda = 1}^{3} \sum_{ \{\vec{k}: |{}^{\lambda}d_{\vec{k}}^{j}| \geq \varepsilon \}} {}^{\lambda}d_{\vec{k}}^{j} \ {}^{\lambda}\bm{\psi}_{\vec{k}}^{j}(\vec{x})
\end{align}
which approximates $f(\vec{x})$ with the spatial error
\begin{align}
\label{eqn:ferror}
    ||f(\vec{x}) - f_{\varepsilon}(\vec{x})||_{\infty} \leq \mathcal{O}(\varepsilon).
\end{align}
Since there exists a one-to-one correspondence between collocation points in the domain and each wavelet basis function, the omission of the wavelet coefficients in the sum corresponds to the omission of collocation points in the computational domain and this procedure results in a sparse, multiresolution spatial discretization with a spatial error bounded by the prescribed tolerance $\varepsilon$.

The wavelet coefficients are mapped back to their representative field values using the matrix operator $\bm{B} = \bm{F}^{-1}$. Both the $\bm{F}$ and $\bm{B}$ matrix operators are sparse, banded, and composed only of the filter coefficients (\emph{i.e.}, $h_{i}, \ \widetilde{h}_{i}, \ g_{i},$ and $\widetilde{g}_{i}$). The use of matrix operators presents the opportunity to replace the cumbersome notation of \cref{eqn:multi} with index notation and implied summation. For instance, the two-dimensional transformations become,
\begin{align}
\label{eqn:FWT_BWT}
	f_{\varepsilon}(\vec{x}) &= \mathbb{d}_{kl} \ \Psi_{kl}(\vec{x}) & \mathrm{with}& & 
	\mathbb{d}_{kl} &= \bm{F}_{kn} \ \bm{F}_{lo} \ f_{no} & \mathrm{and}& &
	 f_{kl} &= \bm{B}_{kn} \ \bm{B}_{lo} \ \mathbb{d}_{no}.
\end{align}
It has been shown that the Deslauriers-Dubuc wavelet family is continuous and differentiable, which allows spatial derivatives to operate directly on the basis functions \cite{Harnish2018AdaptiveControl, Rioul1992SimpleSchemes},
\begin{align}
\label{eqn:D_basis}
    \frac{\partial^{\alpha}}{\partial x_{i}^{\alpha}} f(\vec{x}) &\approx \frac{\partial^{\alpha}}{\partial x_{i}^{\alpha}}\big( \mathbb{d}_{kl} \ \Psi_{kl}(\vec{x}) \big) = \mathbb{d}_{kl} \frac{\partial^{\alpha} \Psi_{kl}(\vec{x})}{\partial x_{i}^{\alpha}}.
\end{align}
As in \cite{Harnish2018AdaptiveControl}, we project the spatial derivative of the basis back onto the wavelet basis and this combination of differentiation and projection transforms \cref{eqn:D_basis} into
\begin{align} 
\label{eqn:D_product}
    \frac{\partial^{\alpha}}{\partial x_{i}^{\alpha}} f(\vec{x}) &\approx \big( \mathcal{D}_{knlo}^{(x_{i}, \alpha)} \mathbb{d}_{no}  \big) \Psi_{kl}(\vec{x}),
\end{align}
where the operator $\mathcal{D}^{(x_{i}, \alpha)}$ is defined in terms of an eigenvector solution and linear combinations of the four types of filter coefficients (\emph{i.e.}, $h_{i}, \ \widetilde{h}_{i}, \ g_{i},$ and $\widetilde{g}_{i}$). Application of the $\mathcal{D}^{(x_{i}, \alpha)}$ operator results in a discrete approximation of the $\alpha^{\mathrm{th}}$ order derivative in the $i$-direction with spatial error
\begin{align}
\label{eqn:D_error}
    \bigg| \bigg| f^{(\alpha)}(\vec{x}) &- \mathcal{D}^{(x_{i}, \alpha)} f_{\varepsilon}(\vec{x}) \bigg| \bigg|_{\infty} \leq \mathcal{O} \left( \varepsilon^{1-\frac{\alpha}{p}} \right).
\end{align}
A derivation of \cref{eqn:D_error} and details on how to assemble the $\mathcal{D}^{(x_{i}, \alpha)}$ operator can be found in our previous publication \cite{Harnish2018AdaptiveControl}. We note that the error estimate in \cref{eqn:D_error} holds only on the dense wavelet grid (\emph{i.e.}, all wavelet coefficients are present) and the sparse wavelet approximation defined in \cref{eqn:multi} approaches this estimate as more of the $\mathcal{D}^{(x_{i}, \alpha)}$ operator's stencil-points are included in the sparse wavelet grid.
\section{Numerical Implementation}
\label{sec:algorithm}
The following algorithm combines the operations defined in \cref{sec:theory} to solve PDEs. To illustrate the process, consider the following unitless model problem,
\begin{align} 
\label{eqn:model}
    \frac{\partial u}{\partial t} +  \vec{V} \cdot \nabla u &= \nu \nabla^{2} u,
\end{align}
with specified initial and boundary conditions, where $\vec{V}$ and $\nu$ are constants. First, the initial condition and operators are projected onto the wavelet basis functions using \cref{eqn:FWT_BWT} and \cref{eqn:D_product}. This process transforms the PDE into an ODE and the two-dimensional semi-discretization of \cref{eqn:model} yields,
\begin{align}
\label{eqn:ode}
    \frac{\mathrm{d} u_{\varepsilon}}{\mathrm{d}t} &= \nu \left[\mathcal{D}^{(x_{1}, 2)} u_{\varepsilon} + \mathcal{D}^{(x_{2}, 2)} u_{\varepsilon} \right] - V_{1} \mathcal{D}^{(x_{1}, 1)} u_{\varepsilon} - V_{2} \mathcal{D}^{(x_{2}, 1)} u_{\varepsilon}.
\end{align}
By retaining collocation points based on the magnitude of the wavelet coefficients, a sparse multiresolution computational grid is created. Furthermore, the spatial accuracy of the initial condition is bounded according to \cref{eqn:ferror}.

Considering that the solution of the PDE may advect and evolve over time, the algorithm must be able to insert and discard points as needed to ensure that the spatial accuracy remains bounded at each time step. Therefore, we use the predictor-corrector strategy described in \cite{Harnish2018AdaptiveControl} to predict where new collocation points will be needed at the next time step and iteratively correct that prediction until the necessary grid is obtained. Specifically, an explicit, embedded, Runge-Kutta time integration scheme \cite{RK23Int1989} is used to progress the solution from the time step $n$ to a trial time step $n+1^{*}$. This converts the ODE in \cref{eqn:ode} into a system of algebraic equations which update $u_{\varepsilon}^{(n)}$ to the trial time step $u_{\varepsilon}^{(n + 1^{*})}$ while providing an estimate of the temporal error and adjusting the time-step size $\Delta t$ such that the temporal error is of the same order as the spatial error (\emph{i.e.}, $\mathcal{O}(\varepsilon)$). At each stage of the Runge-Kutta integration, Dirichlet and Neumann boundary conditions are applied either directly using the procedure developed by Carpenter \emph{et al.} \cite{Carpenter1995BC} (\emph{i.e.}, for Dirichlet data) or by using the penalty formulation \cite{Hesthaven1996BC} that will modify \cref{eqn:ode} on the boundary (\emph{i.e.}, for Dirichlet and/or Neumann data). The magnitudes of the wavelet coefficients at the new time will determine if the grid prediction must be corrected. If so, the trial time step $u_{\varepsilon}^{(n + 1^{*})}$ is discarded and the computational grid at time step $n$ is supplemented with new collocation points and steps are repeated as described in \cite{Harnish2018AdaptiveControl}. When the trial time step is accepted as the true time step, $u_{\varepsilon}^{(n + 1)} = u_{\varepsilon}^{(n + 1^{*})}$, some wavelet coefficients are no longer needed to satisfy the error bounds and are therefore pruned from the sparse computational grid as it evolves with the solution of the PDE \cite{Harnish2018AdaptiveControl}.

This algorithm has been implemented in the Multiresolution Wavelet Toolkit (MRWT) written using modern C++ and is multithreaded using OpenMP. The compressed sparse geometry is stored in a sorted coordinate list (COO) matrix format while collocated field and meta-data are stored in an associated struct-of-array layout. Dynamic grid modification requires a merge and sort operation, but is infrequent and performed in bulk. This grid structure is optimized for slice-based stencil operations versus random access in order to leverage temporal and spatial locality, and is trivially vectorizable for right-hand-side computations. The core $\bm{F}$, $\bm{B}$, and $\mathcal{D}^{(x_{i}, \alpha)}$ operators are stored mostly matrix-free and target the grid's slice-based Application Programming Interface (API). These stencil contractions are trivially parallelizable and scale well.
\section{Numerical Examples}
\label{sec:verification_2D}
\begin{figure}[!htb]
\begin{subfigure}[t]{0.45\textwidth}
    \centering
    \includegraphics[width=\linewidth]{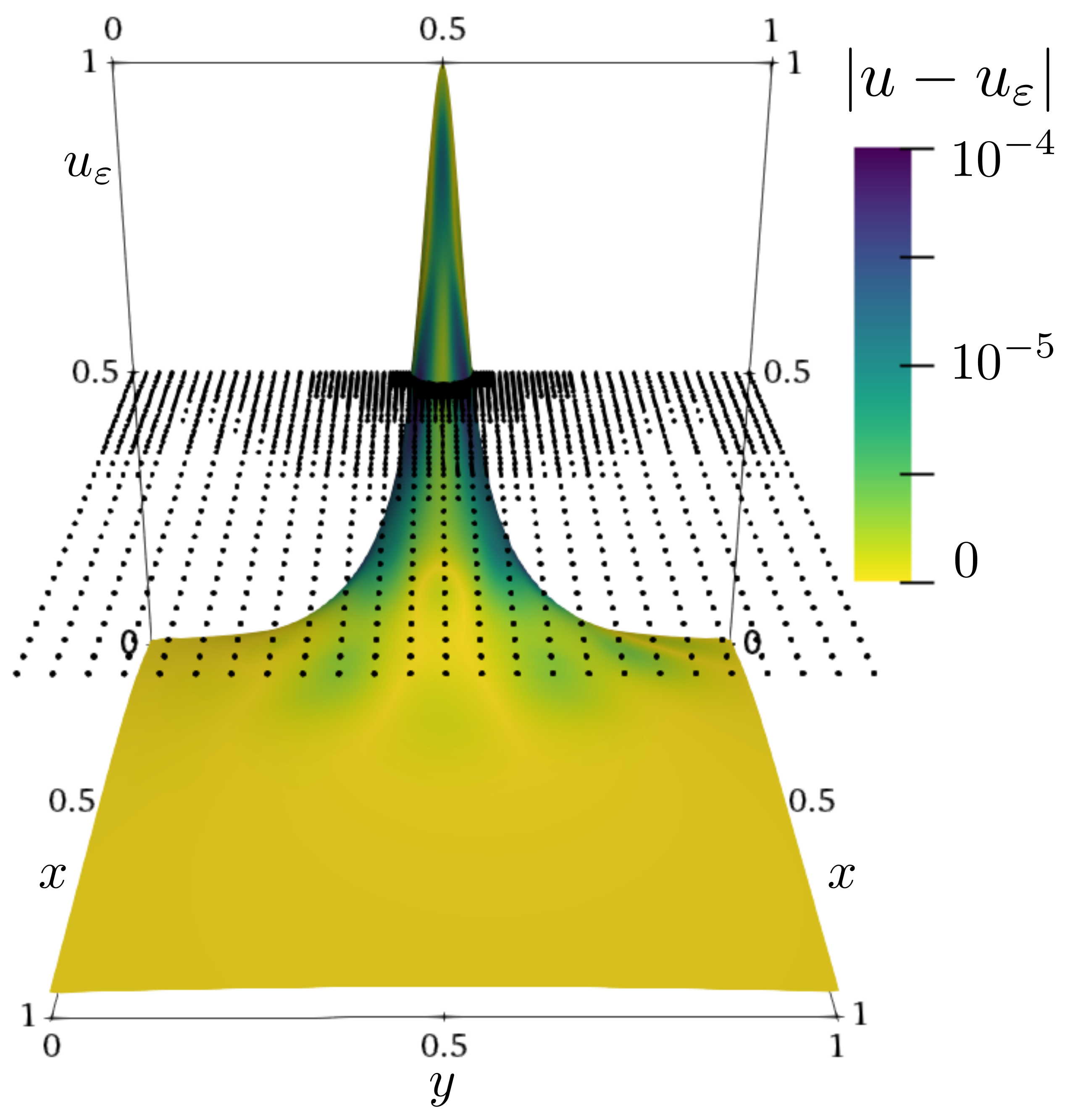}
    \caption{MRWT solution halfway through time using $p = 6$ and $\varepsilon = 10^{-3}$.}
\label{fig:model_solution}
\end{subfigure}
\hfill
\begin{subfigure}[t]{0.48\textwidth}
    \centering
    \includegraphics[width=\linewidth]{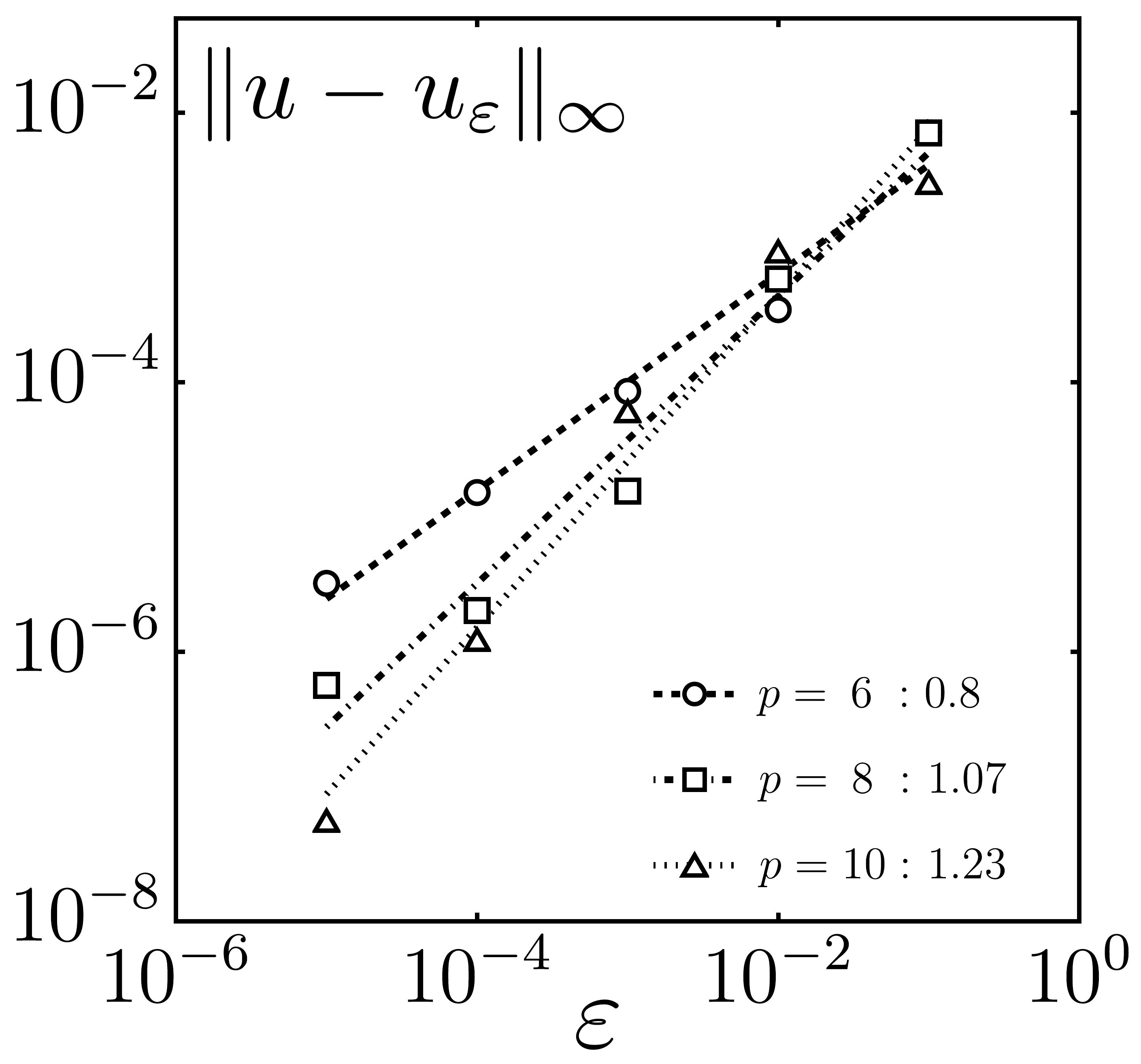}
    \caption{Spatial convergence halfway through time.}
\label{fig:model_convergence}
\end{subfigure}
\centering
\caption{MRWT solution of \cref{eqn:model}, sparse multiresolution grid, and spatial convergence rates obtained halfway through the simulation.}
\label{fig:model_metrics}
\end{figure}

This section provides implementation verification of the algorithm described in \cref{sec:algorithm}. Here, the two-dimensional model problem given by \cref{eqn:model} is solved on the spatial domain, $\Omega = (0,1)^{2}$ and the temporal domain, $T = (0,1/2)$. The specified Dirichlet condition and the initial condition are chosen so that the exact solution is given by,
\begin{align}
\label{eqn:model_exact}
    u(\vec{x}, t) = \frac{5 \nu \left( x_{1} + 5 \nu \right)}{\left( x_{1} + 5 \nu \right)^{2} + \left( t - x_{2}\right)^{2}}.
\end{align}
The constants are set to $V_{1} = 0, \ V_{2} = 1,$ and $\nu = 1/100$. \Cref{eqn:ode} is integrated using the embedded $\mathcal{O}(\Delta t^{2})$ and $\mathcal{O}(\Delta t^{3})$ explicit Runge-Kutta method developed in \cite{RK23Int1989} to maintain $\mathcal{O}(\varepsilon)$ accuracy. \Cref{fig:model_solution} shows the sparse multiresolution grid and the corresponding wavelet approximation, $u_{\varepsilon}$, of the solution halfway through the simulation using wavelet parameters $p = 6$ and $\varepsilon = 10^{-3}$. The exact solution in \cref{eqn:model_exact} allows quantitative error analysis and the color-map in \cref{fig:model_solution} reflects the error $| u - u_{\varepsilon} |$. Convergence rates are calculated by solving \cref{eqn:model} with a variety of wavelet basis $p$ and threshold parameters $\varepsilon$. The convergence of the spatial error is calculated halfway through the simulation. \Cref{fig:model_convergence} shows the error in agreement with estimates \cref{eqn:ferror,eqn:D_error}.

Having demonstrated the mathematical correctness of the numerical method, we now exercise the physics simulation capabilities by solving the coupled system of nonlinear PDEs given by the conservation of mass, momentum, and energy:
\begin{align}
\label{eqn:NS_density}
    \frac{\partial \rho}{\partial t} &= - \nabla \cdot \left( \rho \vec{v}\right), \\
\label{eqn:NS_velocity}
    \frac{\partial}{\partial t} (\rho \vec{v}) &= -\nabla \cdot \left( \rho \vec{v} \otimes \vec{v} - \bm{\sigma}\right) + \rho \vec{b}, \\
\label{eqn:NS_energy}
    \frac{\partial}{\partial t} (\rho \widetilde{e}) &= -\nabla \cdot \left( \rho \widetilde{e} \vec{v} - \bm{\sigma} \vec{v} + \vec{q} \right) + \rho \vec{b} \cdot \vec{v} + \rho r,
\end{align}
where $\widetilde{e} = e + \vec{v} \cdot \vec{v} / 2$. In \cref{eqn:NS_density,eqn:NS_velocity,eqn:NS_energy}, we solve for the density $\rho$, velocity $\vec{v}$, and specific internal energy $e$. This system requires closure equations to describe the Cauchy stress tensor $\bm{\sigma}$, the specific internal energy $e$, and the heat flux $\vec{q}$. In this work, the source terms are set to zero (\emph{i.e.}, $\vec{b} = \vec{0}$ and $r = 0$), we define the stress tensor with the Newtonian fluid constitutive equation, and assume a calorically perfect ideal gas with Fourier’s law of heat conduction.  Additionally, the material parameters are set according to the values in \cref{tab:materialParameters}.
\begin{table}[!htb]
\centering
\begin{tabular}{ |c|c|c| }
    \hline
    Variable & Name & Value \\ 
    \hline
    $\gamma$ & Ratio of specific heats & $7/5$\\
    \hline
    $\mu$ & Dynamic viscosity & $1.9 \times 10^{-5} \ \mathrm{Pa} \cdot \mathrm{s}$ \\ 
    \hline
    $\kappa$ & Thermal conductivity & $2.55 \times 10^{-2} \ \mathrm{W} / (\mathrm{m} \cdot \mathrm{K})$\\ 
    \hline
    $c_{v}$ & Constant volume specific heat & $7.18 \times 10^{2} \ \mathrm{J} / (\mathrm{kg} \cdot \mathrm{K})$ \\ 
    \hline
\end{tabular}
\caption{Material parameters for dry air at room temperature.}
\label{tab:materialParameters}
\end{table}

\begin{figure}[!htb]
\begin{subfigure}[t]{0.49\textwidth}
    \includegraphics[width=0.9\textwidth]{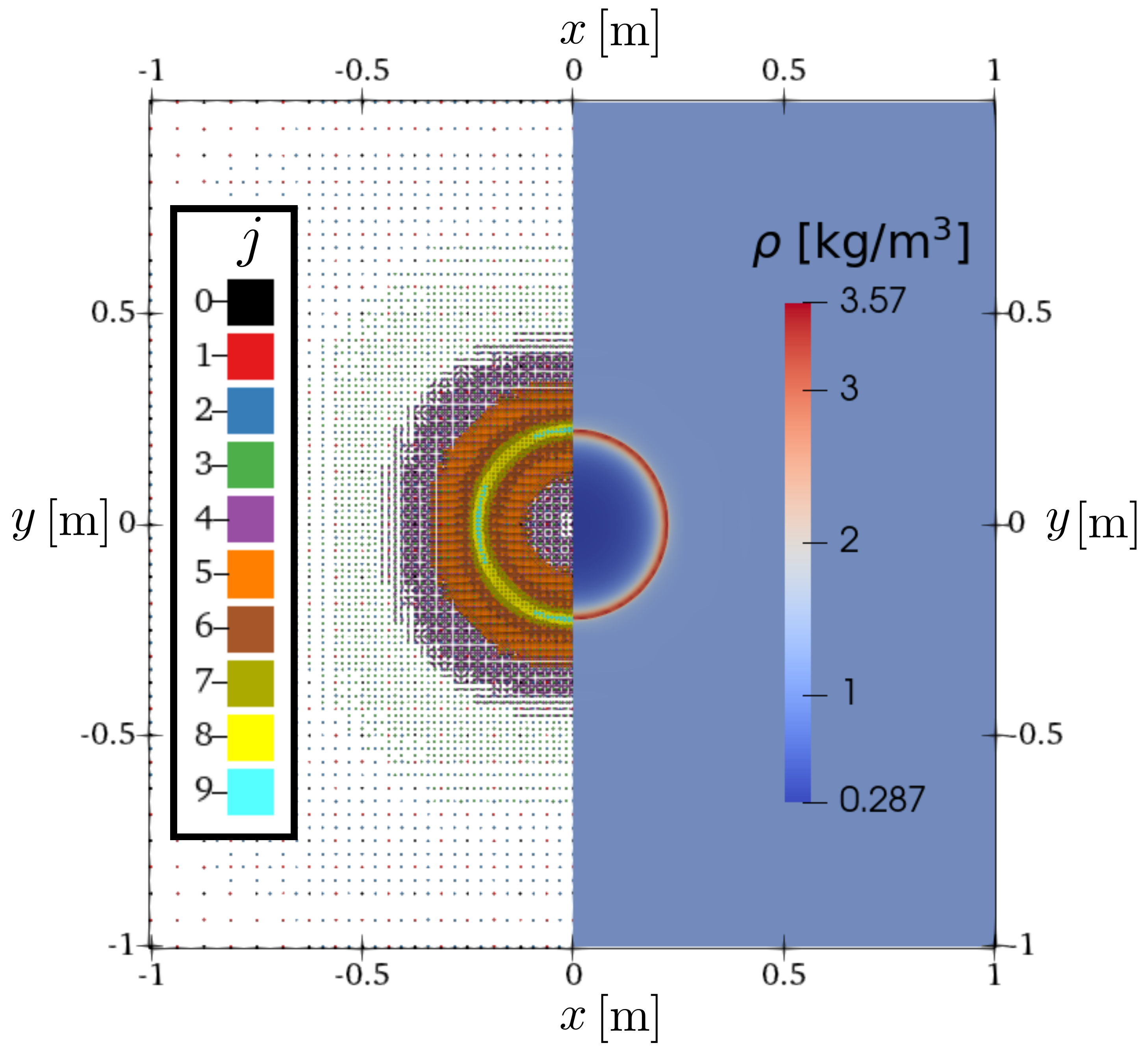}
    \caption{MRWT solution of the density field $\rho$.}
    \label{fig:sedov_grid_density}
\end{subfigure}
\hfill
\begin{subfigure}[t]{0.50\textwidth}
    \includegraphics[width=0.9\textwidth]{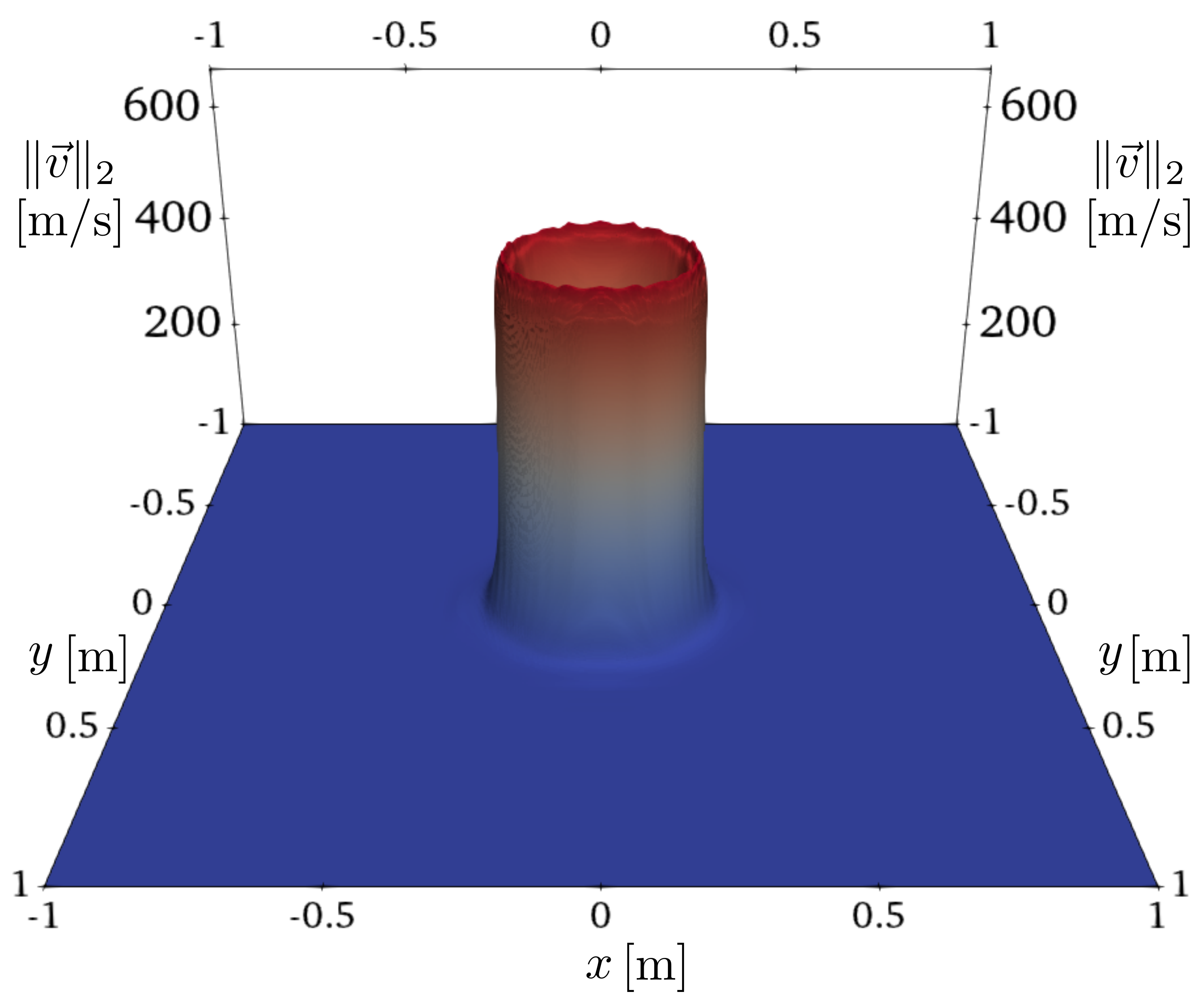}
    \caption{MRWT solution of the velocity field $\| \vec{v} \|_{2}$.}
    \label{fig:sedov_velocity}
\end{subfigure}
\centering
\caption{Sparse multiresolution grid and numerical solution at $t = 133.902 \ \mu$s obtained using $p = 8$ and $\varepsilon = 10^{-2}$. In (a), the grid points are colored according to their resolution level $j$. The reader is referred to the online version of this article for clarity regarding the color in this figure. In (b), the maximum velocity is approximately $568$ m/s.
}
\label{fig:sedov}
\end{figure}

This model has been used to describe the evolution of a Taylor-Sedov blast wave \cite{Paolucci2014WAMR:Algorithm}, where energy is deposited in a compressible fluid leading to the development of a spherical shock wave. The initial condition is made continuous by way of a Gaussian profile for the initial pressure, with an overpressure peak of $2$ MPa and a standard deviation of $1/(10\sqrt{2})$ m. The semi-discretized \cref{eqn:NS_density,eqn:NS_velocity,eqn:NS_energy} are integrated using the embedded $\mathcal{O}(\Delta t^{4})$ and $\mathcal{O}(\Delta t^{5})$ explicit Runge-Kutta method developed in \cite{RKF1969}. The temporal discretization, $\Delta t$, is chosen adaptively to retain $\mathcal{O}(\varepsilon)$ accuracy. The boundary conditions are set to maintain the initial conditions and the simulation is stopped before the developing shock wave interacts with the computational boundary. \Cref{fig:sedov} shows the numerical solutions to \cref{eqn:NS_density,eqn:NS_velocity,eqn:NS_energy} at time $t = 133.902 \ \mu$s, generated with wavelet parameters $p = 8$ and $\varepsilon = 10^{-2}$.

For the numerical solution in \cref{fig:sedov}, the MRWT discretization of the initial condition required only two resolution levels (\emph{i.e.}, $j = 2$) which resulted in $31.250$ mm between the closest collocation points at time $t = 0$. As the internal energy converted into kinetic energy, MRWT automatically refined the grid near regions of the developing shock wave. As shown in \cref{fig:sedov_grid_density}, MRWT predicted nine resolution levels (\emph{i.e.}, $j = 9$) at time $t = 133.902 \ \mu$s which resulted in $0.244$ mm between the closest collocation points. A dense discretization at this length scale would require over 67 million collocation points whereas the MRWT solution in \cref{fig:sedov} only contains 312,793 collocation points, resulting in a compression ratio greater than 200. Moreover, the sparse multiresolution spatial discretization maintains symmetry and adapts to follow features as they evolve through the domain.

\section{Conclusions}
Our proposed wavelet based algorithm has the potential to improve computational science and engineering applications requiring resolution across multiple spatial and temporal scales. In particular, our method is well suited for situations where fine scale features are dynamic. Whereas traditional numerical methods require costly remeshing/refining procedures, our method leverages the properties of wavelet basis functions to automatically adapt the computational domain as needed to accurately resolve features.

We have demonstrated that our implementation, MRWT, is capable of solving multidimensional PDEs with error controlled by the threshold parameter $\varepsilon$. This work advances the state of wavelet based algorithms by exploiting the regularity of the Deslauriers-Dubuc family of wavelets to evaluate spatial derivatives directly on the wavelet basis functions. Additionally, we have verified that our predictor-corrector procedure is able to solve initial-boundary value problems using compressed data on sparse multiresolution discretizations in finite domains. Furthermore, we have provided error estimates for each wavelet operation and have shown that our numerical solutions have bounded error at each time step, with convergence rates in agreement with the theoretical estimates.


\section*{Acknowledgment}
This work was supported by Los Alamos National Laboratory (LANL) under award numbers $369229$ \& $370985$.


\begin{thebibliography}{102}
\bibitem{ocean}
Ringler T, Petersen M, Higdon R, Jacobsen D, Jones P, Maltrud M 2013 {\it Ocean Modelling} {\bf 69} 211--232

\bibitem{detonation}
Cai X, Liang J, Deiterding R, Che Y, Lin Z 2016 {\it International Journal of Hydrogen Energy} {\bf 41} 3222--3239

\bibitem{asteroid}
Boslough M, Jennings B, Carvey B, Fogleman W 2015 {\it Procedia Engineering} {\bf 103} 43--51

\bibitem{matousJCP}
Matouš K, Geers M, Kouznetsova V, Gillman A 2017 {\it J. Comput. Phys.} {\bf 330} 192--220

\bibitem{supernova}
Malone C et al. 2014 {\it The Astrophysical Journal} {\bf 782}(1) 1--24

\bibitem{Berger1984AdaptiveEquations}
Berger M, Oliger J 1984 {\it J. Comput. Phys.} {\bf 53}(3) 484--512

\bibitem{Fatkullin2001AdaptiveProblems}
Fatkullin I, Hesthaven J 2001 {\it Journal of Scientific Computing} {\bf 16}(1) 47--67

\bibitem{Brandt1977Multi-LevelProblems}
Brandt A 1977 {\it Mathematics of Computation} {\bf 31}(138) 333

\bibitem{Hackbusch1978OnEquations}
Hackbusch W 1978 {\it Computing} {\bf 20} 291--306

\bibitem{matousDewen}
Yushu D, Matouš M 2020 {\it J. Comput. Phys.} {\bf 406} 109165

\bibitem{chimera}
Benek J et al. 1986 {\it Arnold Engineering Development Center Technical Report} {\bf 85-64}

\bibitem{Dong2003P-refinementP-threads}
Dong S, Karniadakis G 2003 {\it CMAME} {\bf 192}(19) 2191--2201
    
\bibitem{Gui1986TheDimension}
Gui W, Babu{\v{s}}ka I 1986 {\it Numerische Mathematik} {\bf 49} 613--657
    
\bibitem{Gui1986TheDimensionb}
Gui W, Babu{\v{s}}ka I 1986 {\it Numerische Mathematik} {\bf 49} 577--612
    
\bibitem{Rajagopal2007AInterfaces}
Rajagopal A, Sivakumar S 2007 {\it Computational Mechanics} {\bf 41} 49--72

\bibitem{Jawerth1994AnAnalyses}
Jawerth B, Sweldens W 1994 {\it SIAM Review} {\bf 36}(3) 377--412

\bibitem{Schneider2010WaveletDynamics}
Schneider K, Vasilyev O 2010 {\it Annual Review of Fluid Mechanics} {\bf 42}(1) 473--503

\bibitem{Liandrat1990ResolutionApproximation}
Liandrat J, Tchamitchian Ph 1990 {\it NASA ICASE Report 90-83}
    
\bibitem{Beylkin1997OnBases}
Beylkin G, Keiser J 1997 {\it J. Comput. Phys.} {\bf 132}(2) 233--259

\bibitem{Bertoluzza1996AdaptiveEquation}
Bertoluzza S 1996 {\it Transport Theory and Statistical Physics} {\bf 25} (3-5) 339--352
    
\bibitem{Ueno2003AProperty}
Ueno T,  Ide T, Okada M 2003 {\it Bulletin des Sciences Mathematiques} {\bf 127}(6) 569--583
    
\bibitem{Qian1993WaveletsEquations}
Qian S, Weiss J 1993 {\it J. Comput. Phys.} {\bf 106} 155--175
    
\bibitem{Kong2016NonlinearTechnique}
Kong F, Kougioumtzoglou I, Spanos P, Li S 2016 {\it Int. J. Multiscale Comput. Eng.} {\bf 14}(3) 255--272
    
\bibitem{RodyEtAl}
van Tuijl R, Harnish C, Matou{\v{s}} K, Remmers J, Geers M 2018 {\it Computational Mechanics} {\bf 63} 535--554
    
\bibitem{Paolucci2014WAMR:Algorithm}
Paolucci S, Zikoski Z, Wirasaet D 2014 {\it J. Comput. Phys.} {\bf 272} 814--841
    
\bibitem{Paolucci2014WAMR:Algorithmb}
Paolucci S, Zikoski Z, Grenga T 2014 {\it J. Comput. Phys.} {\bf 272} 842--864
    
\bibitem{Nejadmalayeri2015ParallelPDEs}
Nejadmalayeri A, Vezolainen A,  Brown-Dymkoski E, Vasilyev O 2015 {\it J. Comput. Phys.} {\bf 298} 237--253
    
\bibitem{Dubos2013AGrids}
Dubos T, Kevlahan N 2013 {\it Quarterly Journal of the Royal Meteorological Society} {\bf 139}(677) 1997--2020
    
\bibitem{Sakurai2017CoherentWavelets}
Sakurai T {\it et al.} 2017 {\it Journal of Turbulence} {\bf 18}(4) 352--372
    
\bibitem{Frohlich1994AnComputations}
Fr{\"{o}}hlich J, Schneider K 1994 {\it European Journal of Mechanics} B {\bf 13}(4) 439--471
    
\bibitem{Goedecker2009WaveletsPhysics}
Goedecker S 1998 (Lausanne : Presses Polytechniques et Universitaires Romandes)
    
\bibitem{Iqbal2014AnModeling}
Iqbal A, Jeoti V 2014 {\it Radio Engineering} {\bf 23}(4) 987--996
    
\bibitem{Le2015Reduced-orderWinds}
Le T, Caracoglia L 2015 {\it Journal of Sound and Vibration} {\bf 344} 179--208
    
\bibitem{Lin2001ConnectionEquation}
Lin E, Zhou X 2001 {\it Journal of Computational and Applied Mathematics} {\bf 135}(1) 63--78
    
\bibitem{Holmstrom1999SolvingWavelets}
Holmstr{\"{o}}m M 1999 {\it SIAM J. Sci. Comput.} {\bf 21}(2) 405--420

\bibitem{Harnish2018AdaptiveControl}
Harnish C, Matou{\v{s}} K, Livescu D 2018 {\it Int. J. Multiscale Comput. Eng.} {\bf 16}(1) 19--43

\bibitem{Daubechies1992TenDaubechies}
Daubechies I 1992 {\it Ten Lectures on Wavelets} (Philadelphia : Society for Industrial and Applied Mathematics)
    
\bibitem{Cohen2000MultiscaleDomains}
Cohen A, Dahmen W, DeVore R 2000 {\it Transactions of the American Mathematical Society} {\bf 352}(8) 3651--3685
    
\bibitem{Bacry1992AEquations}
Bacry E, Mallat S, Papanicolaou G 1992 {\it Mod{\'{e}}lisation math{\'{e}}matique et analyse num{\'{e}}rique} 793--834
    
\bibitem{deVilliers2003Dubuc--DeslauriersInterval}
de Villiers J, Goosen K, Herbst B 2003 {\it SIAM Journal on Mathematical Analysis} {\bf 35}(2) 423--452
    
\bibitem{Rioul1992SimpleSchemes}
Rioul O 1992 {\it SIAM Journal on Mathematical Analysis} {\bf 23}(6) 1544--1576
    
\bibitem{RK23Int1989}
Bogacki P, Shampine L 1989 {\it Applied Mathematics Letters} {\bf 2}(4) 321--325
    
\bibitem{Carpenter1995BC}
Carpenter M, Gottlieb D, Abarbanel S, Don W 1995 {\it SIAM J. Sci. Comput.} {\bf 16}(6) 1241--1252
    
\bibitem{Hesthaven1996BC}
Hesthaven J, Gottlieb D 1996 {\it SIAM J. Sci. Comput.} {\bf 17}(3) 579--612

\bibitem{RKF1969}
Fehlberg E,  1969 {\it NASA Technical Report} {\bf R-315}

\end{thebibliography}
\end{document}